\documentclass{article}
\usepackage{amsmath, amsthm, amssymb}
\usepackage{makeidx}
\usepackage{amsbsy}
\usepackage{graphicx}
\begin{document}
\def\theequation{\arabic{section}.\arabic{equation}}
\setcounter{page}{1}

\newcommand{\C}{C^{-1}}
\newcommand{\bt}{\tilde{\boldsymbol{\beta}}}
\newcommand{\bh}{\hat{\boldsymbol{\beta}}}
\newcommand{\V}{\mathbf{V}}
\newcommand{\X}{\mathbf{X}}
\newcommand{\bV}{\boldsymbol{V}}

\newcommand{\by}{\boldsymbol{y}}
\newcommand{\ba}{\boldsymbol{a}}
\newcommand{\bz}{\boldsymbol{z}}
\newcommand{\bY}{\boldsymbol{Y}}
\newcommand{\bA}{\boldsymbol{A}}
\newcommand{\bx}{\boldsymbol{x}}
\newcommand{\bw}{\boldsymbol{w}}
\newcommand{\bX}{\boldsymbol{X}}
\newcommand{\bB}{\boldsymbol{B}}
\newcommand{\bC}{\boldsymbol{C}}
\newcommand{\bH}{\boldsymbol{H}}
\newcommand{\bM}{\boldsymbol{M}}
\newcommand{\bQ}{\boldsymbol{Q}}
\newcommand{\bbh}{\boldsymbol{h}}
\newcommand{\bR}{\boldsymbol{R}}
\newcommand{\bS}{\boldsymbol{S}}
\newcommand{\bI}{\boldsymbol{I}}
\newcommand{\bW}{\boldsymbol{W}}
\newcommand{\bSi}{\boldsymbol{\Sigma}}
\newcommand{\0}{\boldsymbol{0}}

\newcommand{\bb}{\boldsymbol{\beta}}
\newcommand{\be}{\boldsymbol{\epsilon}}
\newcommand{\bet}{\boldsymbol{\eta}}
\newcommand{\bd}{\boldsymbol{\delta}}
\newcommand{\bmu}{\boldsymbol{\mu}}
\newcommand{\W}{\mathcal{W}}
\newcommand{\A}{\mathcal{A}}
\newcommand{\B}{\mathcal{B}}
\newcommand{\E}{\mathcal{E}}
\newcommand{\N}{\mathcal{N}}
\newcommand{\D}{\mathcal{\boldsymbol{D}}}
\newcommand{\beq}{\begin{eqnarray}}
\newcommand{\eeq}{\end{eqnarray}}
\newtheorem{thm}{Theorem}[section]
\newtheorem{cor}{Corollary}[section]
\newtheorem{lem}{Lemma}[section]
\newtheorem{rem}{Remark}[section]

\newcommand{\nl}{\mathcal{L}_n}

\title{\bf Regression Model With Elliptically Contoured Errors}

\vspace{3.5cm}

\author{{ M. Arashi\footnote{Corresponding Author ({\em Email: m\_arashi\_stat@yahoo.com})}  $^1$, A. K. Md E. Saleh
$^2$ and S. M. M. Tabatabaey$^3$}\vspace{.5cm}
\\\it
$^{1}$Faculty of Mathematics, Shahrood University of Technology,
\\\vspace{.2cm}\it P.O. Box 316-3619995161, Shahrood, Iran.
\\\it $^2$School of Mathematics and Statistics,
Carleton University, Ottawa\\\it $^{3}$Department of Statistics,
School of Mathematical Sciences\\\vspace{.5cm} \it Ferdowsi
University of Mashhad, Mashhad, Iran}

\date{}
\maketitle

\begin{quotation}
\noindent {\it Abstract:} For the regression model
$\by=\bX\bb+\be$ where the errors follow the elliptically
contoured distribution (ECD), we consider the least squares (LS),
restricted LS (RLS), preliminary test (PT), Stein-type shrinkage
(S) and positive-rule shrinkage (PRS) estimators for the
regression parameters, $\bb$.\\ We compare the quadratic risks of
the estimators to determine the relative dominance properties of
the five estimators.
\par

\vspace{9pt} \noindent {\it Key words and phrases:} Elliptically
contoured distribution, Multivariate Student's t, Positive-rule
shrinkage estimator, Preliminary test estimator, Signed measure,
Stein-type shrinkage estimator.
\par
\end{quotation}\par

\section{Introduction}
The most important model belonging to the class of general linear
hypotheses is the {\it multiple regression} model. The general
purpose of {\it multiple regression} is to learn more about the
relationship between several independent or predictor variables
and a dependent or criterion variable.

Consider the multiple regression model
\begin{equation}\label{eq11}
\by=\bX\bb +\be
\end{equation}
where $\by$ is an $n$-vector of responses, $\bX$ is an $n\times p$
non-stochastic design matrix with full column rank $p$,
$\bb=(\beta_1,\cdots,\beta_p)'$ is $p$-vector of regression
coefficients and $\be=(\epsilon_1,\cdots,\epsilon_n)'$ is the
$n$-vector of random errors distributed according to the law
belonging to the class of elliptically contoured distributions
(ECDs), $\E_n(\0,\sigma^2\bV,\psi)$ for $\sigma\in\mathbb{R}^+$
and un-structured known matrix $\bV\in S(n)$, where $S(n)$ denotes
the set of all positive definite matrices of order $(n\times n)$
with the following characteristic function
\begin{equation}\label{eq12}
\phi_{\be}(\boldsymbol{t}) = \psi\left(\sigma^2\boldsymbol{t}'\bV
\boldsymbol{t}\right)
\end{equation}
for some functions $\psi: [0,\infty)\rightarrow\mathbb{R}$ say
characteristic generator (Fang et al., 1990).

If $\be$ has a density, then it is of the form
\begin{eqnarray}\label{eq13}
f(e)& \propto
&|\sigma^2\mathbf{V}|^{-\frac{1}{2}}g\left(\frac{1}{\sigma^2}\;\be'\bV^{-1}\be\right)
\end{eqnarray}
where $g(.)$ is a non-negative function over $\mathbb{R}^+$ such
that $f(.)$ is a density function w.r.t (with respect to) a
$\sigma$-finite measure $\mu$ on $\mathbb{R}^p$. In this case,
notation $\be\sim\E_{n}(\0,\sigma^2\bV,g)$ would probably be used.

It is sometimes difficult to have complete analysis of the
regression model with ECD errors of the type \eqref{eq12} or
\eqref{eq13}. To overcome such difficulties, one may consider any
of the three sub-classes of ECDs, namely,
\begin{enumerate}
\item[(i)] Scale mixture of normal distributions,
\item[(ii)] Laplace class of mixture of normal distributions, and
\item[(iii)] Signed measure mixture of normal distributions.
\end{enumerate}
General formula for the above mixture of distributions is given by
\begin{eqnarray}\label{eq14}
f_{\be}(\bx)=\int_0^\infty
\W(t)\phi_{\N_n(\0,t^{-1}\sigma^2\bV)}(\bx)dt,
\end{eqnarray}
where $\phi_{\N_n(\0,t^{-1}\sigma^2\bV)}(.)$ is the pdf
(probability density function) of $\N_n(\0,t^{-1}\sigma^2\bV)$.\\
\begin{enumerate}
\item[(a)] If
\begin{equation}\label{eq15}
\W(\tau)=2\left(\Gamma(\gamma/2)\right)^{-1}
\left(\frac{\gamma\sigma^2}{2}\right)^{\gamma/2}\tau^{-(\gamma+1)}
e^{-\frac{\gamma\sigma^2}{2\tau^2}},\qquad
0<\gamma,\sigma^2,\tau<\infty
\end{equation}
then we have
\begin{equation}\label{eq16}
f(\be)=\frac{\Gamma\left(\frac{n+\gamma}{2}\right)|\bV|^{-\frac{1}{2}}}{(\pi\gamma)^{n/2}\Gamma\left(\gamma/2\right)\sigma^n}
\left(1+\frac{\be'\bV^{-1}\be}{\gamma\sigma^2}\right)^{-\frac{1}{2}(n+\gamma)},
\end{equation}
where $E(\be)=\mathbf{0}$ and
$E(\be\be')=\frac{n\gamma\sigma^2}{\gamma-2}\bV=\sigma_e^2$ for
$\gamma>2$.
\item[(b)] Chu (1973) considered
\begin{eqnarray}\label{eq17}
\W(t)=(2\pi)^{\frac{n}{2}}|\sigma^2\bV|^{\frac{1}{2}}t^{-\frac{p}{2}}\mathcal{L}^{-1}[f(s)],
\end{eqnarray}
$\mathcal{L}^{-1}[f(s)]$ denotes the inverse Laplace transform of
$f(s)$ with $s=[\bx'(\sigma^2\bV)^{-1}\bx/2]$. For some examples
of $f(.)$ and $\W(.)$ see Arashi et al. (2010).

The inverse Laplace transform of $f(.)$ exists provided that the
following conditions are satisfied.
\begin{enumerate}
\item[(i)]$f(t)$ is differentiable when $t$ is sufficiently large.
\item[(ii)]$f(t)=o(t^{-m})$ as $t\rightarrow\infty$, $m>1$.
\end{enumerate}
Although, it is rather difficult to derive the inverse Laplace
transform of some functions, we are able to handle it for many
density generators of elliptical densities. We refer the readers
to Debnath and Batta (2007) for more specific details.

The mean of $\be$ is the zero-vector and the covariance-matrix of
$\be$ is
\begin{eqnarray}\label{eq18}
\bSi_{\be}=Cov(\be)&=&\int_0^\infty Cov(\be|t)\W(t)dt\nonumber\\
&=&\int_0^\infty \W(t)Cov\left\{\N_p(\0,t^{-1}\sigma^2\bV)\right\}dt\nonumber\\
&=&\left(\int_0^\infty t^{-1}\W(t)dt\right)\sigma^2\bV,
\end{eqnarray}
provided the above integral exists.\\
Comparing the models \eqref{eq13} and \eqref{eq14}, since
$\bSi_{\be}=Cov(\be)=-2\psi'(0)\sigma^2\bV$, using \eqref{eq16} we
can conclude that
\begin{eqnarray*}
-2\psi'(0)=\int_0^\infty t^{-1}\W(t)dt.
\end{eqnarray*}

Now suppose that $\bX\sim \E_n(\bmu,\bV,g)$. Then it is important
to point out that since $\int_{\bx} f(\bx)d\bx=1$, using Fubini's
theorem we have
\begin{eqnarray*}
1&=&\int_{\bx} \int_0^\infty \W(t)\phi_{\N_n(\bmu,t^{-1}\bV)}(\bx)dt d\bx\\
&=& \int_0^\infty \W(t)\int_{\bx}\phi_{\N_n(\bmu,t^{-1}\bV)}(\bx)d\bx dt\\
&=&\int_0^\infty \W(t)dt.
\end{eqnarray*}
Thus for nonnegative function $\W(.)$, it is a density. For
nonnegative function $\W(.)$, the elliptical models can be
interpreted as a scale mixture of normal distributions.
\item[(c)] Srivastava and Bilodeau (1989) considered the signed measure, $\W(t)$ such
that
\begin{eqnarray}\label{eq19}
(i)&&\int_0^\infty t^{-1}\W^+(dt)<\infty,\nonumber\\
(ii)&&\int_0^\infty t^{-1}\W^-(dt)<\infty,
\end{eqnarray}
where $\W^+ - \W^-$ is the Jordan decomposition of $\W$ in
positive and negative parts. Note that from $(i)-(ii)$ of
\eqref{eq19},
\begin{eqnarray}\label{eq110}
\int_0^\infty t^{-1}\W(dt)<\infty
\end{eqnarray}
and thus, $Cov(\be)$ exists under the sub-class defined above.\\
This subclass contains the subclass defined by (b).
\end{enumerate}
\begin{rem}
Regarding the above classifications, we should take the following
notes:
\begin{enumerate}
\item[1.] In all the above classes we have
\begin{eqnarray*}
\bSi_{\be}=-2\psi'(0)\sigma^2\bV=\left(\int_0^\infty
t^{-1}\W(t)dt\right)\sigma^2\bV,
\end{eqnarray*}
resulting in $-2\psi'(0)=\int_0^\infty t^{-1}\W(t)dt$.
\item[2.] The subclass (a) is neither contained in subclass (b) nor in the
subclass (c). However, subclass (b) in contained in the
subclass(c). Thus, all the implications about the subclass (c) can
be used for the subclass (b).
\item[3.] For the subclass (c) we can assure that $-2\psi'(0)=\int_0^\infty
t^{-1}\W(t)dt$ exists. However it may not exist for the subclass
(b).
\end{enumerate}
\end{rem}

Some of the well-known members of the class of ECDs are the
multivariate normal, Kotz Type, Pearson Type II \& VII,
multivariate Student's t, multivariate Cauchy, Logistic, Bessel
and generalized slash distributions. Dating back to Kelker (1970),
there are many known results concerning ECDs, in particular the
mathematical properties and its application to statistical
inference. These results have been put forward by Cambanis et al.
(1981), Muirhead (1982), Fang et al. (1990) and Gupta and Varga
(1993) among others .

The object of the paper is the estimation of the regression
parameters, $\bb=(\beta_1,\cdots,\beta_p)'$ when it is suspected
that $\bb$ may belong to the sub-space defined by $\bH\bb=\bbh$
where $\bH$ is a $q\; \mathsf{x}\; p$ matrix of constants and
$\bbh$ is a q-vector of known constants with focus on the
Stein-type estimators of $\bb$ in addition to preliminary test
estimator (PTE) based on the error distributions belonging to the
subclass(c) which includes subclass(b) described earlier.

The recent book of Saleh (2006) dealing with the proposed
estimators (chapter 7) presents an overview on the topic under
normal as well as nonparametric theory covering many standard
statistical models. Tabatabaey (1995) and Arashi and Tabatabaey
(2008) covers the theory with spherically symmetric distribution
of errors developing many interesting calculations. For some
systematic work have been done so far in the context of Stein-type
estimations see Srivastava and Bilodeau (1989) and Arashi and
Tabatabaey (2009) under elliptical symmetry.

We organize our paper as follows: Section 2, contains the
estimation and the test of hypothesis along with proposed
estimators of $\bb$, section 3 deals with the bias, risk and MSE
expressions of the proposed estimators while the analysis of the
risks and comparisons are presented in section 4. Concluding
remarks are presented in section 5.

\setcounter{equation}{0}
\section{Estimation and Test of Hypothesis}
In this section, we present the estimate of $\bb$ and $\sigma^2$
under least square (LS) theory. Further we discuss the problem of
testing the general linear hypothesis, $\bH\bb=\bbh$. The test of
this hypothesis covers many special
cases considered in practical situations.\\
Using standard conditions, it is well-known that the generalized LS (GLS) estimator
of $\bb$ is
\begin{equation}\label{eq21}
\tilde{\bb}=(\bX'\bV^{-1}\bX)^{-1}\bX'\bV^{-1}\by=\bC^{-1}\bX'\by,\qquad
\bC=\bX'\bV^{-1}\bX.
\end{equation}
Under elliptical assumptions, its distribution is
$\E_p(\bb,\sigma^2\bC^{-1},g)$. Similarly the estimate of the
$\sigma^2$ is
\begin{equation}\label{eq22}
\tilde{\sigma}^2=\frac{1}{n}(\by-\bX\tilde{\bb})'\bV^{-1}(\by-\bX\tilde{\bb}).
\end{equation}
It is easy to show that
\begin{equation}\label{eq23}
S^2=\frac{(\by-\bX\tilde{\bb})'\bV^{-1}(\by-\bX\tilde{\bb})}{n-p}
\end{equation}
is an unbiased estimator of $\sigma^2_{\be}=-2\psi'(0)\sigma^2$.\\
For test of $\bH\bb=\bbh$ (where $q<p$), we first consider the
restricted estimator given by
\begin{equation}\label{eq24}
\hat{\bb}=\bt- \bC \bH'\bV_1(\bH\bt-\bbh),\quad\bV_1=(\bH\bC
\bH')^{-1}.
\end{equation}
It can be directly verified that $\hat{\bb}\sim \E_p\left(\bb-\bC
\bH'\bV_1(\bH\bb-\bbh),\sigma^2\bV_2,g\right)$ for
$\bV_2=\bC(\bI_p-\bH'\bV_1\bH\bC)$. Therefore, we get
\begin{eqnarray}\label{eq25}
E(\bh -\bb)&=&-\bC \bH'\bV_1(\bH\bb-\bbh)\nonumber\\
&=&-\bd\quad(\textrm{say}),\nonumber\\
E(\bh-\bb)'(\bh-\bb)&=&-2\sigma^2\Psi'(0)\boldsymbol{tr}(\bV_2)+\bd'\bd\nonumber\\
&=&\sigma_{\be}^2\boldsymbol{tr}(\bV_2)+\bd'\bd.
\end{eqnarray}
Similarly, under $\bH\bb=\bbh$, the following estimator is
unbiased for $\sigma^2_{\be}$.
\begin{equation}
S^{*2}=\frac{(\by-\bX\bh)'\bV^{-1}(\by-\bX\bh)}{n-p+q},
\end{equation}
from least square's theory.\\
Now we consider the linear hypothesis $\bH\bb=\bbh$ and obtain the
likelihood ratio test (LRT) statistic for the null hypothesis
$H_0:\bH\bb=\bbh$ as well as its non null distribution. The LRT
statistic is a consequence of the results of Anderson et al.
(1986).
\begin{thm}
Let
\begin{eqnarray*}
w&=&\{\bb:\bb \in \mathbb{R}^p,\bH\bb=\bbh,\sigma >0, \bV\in
S(n)\}, \mbox{and}\\
\Omega&=&\{\bb:\bb \in \mathbb{R}^p,\sigma >0, \bV\in S(n)\}
\end{eqnarray*}
Moreover, suppose $y^{\frac{n}{2}}g(y)$ has a finite positive
maximum $y_{g}$. Then under the assumptions of model (1.1) the LRT
for testing the null-hypothesis $H_0:\bH\bb=\bbh$ is given by
\begin{equation}
\mathcal{L}_n=\frac{(\bH\bt-\bbh)'\bV_1(\bH\bt-\bbh)}{qS^2},
\end{equation}
and has the following generalized non-central F-distribution with
pdf
\begin{equation}
\boldsymbol{g}_{q,m}^*(\nl)= \sum_{r\geq 0}
\frac{\left(\frac{q}{m}\right)^{\frac{1}{2}(q+2r)}\nl^{\frac{1}{2}(q+2r-2)}K^{(0)}_r(\Delta_*^2)}
{r!\;
\boldsymbol{B}\left(\frac{q+2r}{2},\frac{m}{2}\right)\left(1+\frac{q}{m}\nl\right)^{\frac{1}{2}(q+m+2r)}},
\end{equation}
where $m=n-p$, $\Delta_*^2=\theta/\sigma_{\be}^2$ for
$\theta=(\bH\bb-\bbh)'\bV_1(\bH\bb-\bbh)$, and the mixing
distribution becomes
\begin{eqnarray}
K^{(0)}_r(\Delta_*^2)&=&[-2\psi'(0)]^r\left(\frac{\Delta_*^2}{2}\right)^r\int^{\infty}_0\frac{t^r}{r!}\;e^{\frac{-t\Delta_*^2[-2\psi'(0)]}{2}}\W(t)dt.
\end{eqnarray}
\end{thm}
For the proof see the Appendix.
\begin{cor}
Under $H_0$, the pdf of $\nl$ is given by
\begin{equation}
\boldsymbol{g}_{q,m}^*(\nl)=
\frac{\left(\frac{q}{m}\right)^{\frac{q}{2}}\nl^{\frac{q}{2}-1}}
{\boldsymbol{B}\left(\frac{q}{2},\frac{m}{2}\right)\left(1+\frac{q}{m}\nl\right)^{\frac{1}{2}(q+m)}},
\end{equation}
which is the central F-distribution with $(q,m)$ degrees of
freedom. For $W(\tau)$ given by \eqref{eq15}, we get the results
produced by Tabatabaey (1995).
\end{cor}
Now, consider the calculations of the probability of that
$\nl\leq F_\alpha$, which gives the power function of the test as
\begin{equation}\label{G}
\boldsymbol{G}_{q,m}^*(F_\alpha;\Delta_*^2)= \sum_{r\geq 0}
\frac{1}{r!}K^{(0)}_r(\Delta_*^2) I_x
\left[\frac{1}{2}(q+2r),\frac{m}{2}\right],
\end{equation}
where $x=\frac{qF_\alpha}{(m+qF_\alpha)}$ and $I_x(a,b)$ is the
incomplete Beta-function,
\begin{equation}
I_x(a,b)=\frac{1}{\boldsymbol{B}(a,b)}\int_0^x u^{a-1}(1-u)^{b-1}
du.
\end{equation}
The function (\ref{G}) stands for the power function at
$\alpha$-level of significance and may be called the generalized
non-central F-distribution cdf (cumulative distribution function)
of the statistic $\nl$.

Similarly, the cdf of a generalized non-central chi-square
distribution with $\gamma$ d.f. may be written as
\begin{equation}\label{H}
\mathcal{H}_\gamma^*(x;\delta^2)=\sum_{r\geq 0}
\frac{1}{r!}K^{(0)}_r(\Delta_*^2)\mathcal{H}_{\gamma+2r}(x;0),
\end{equation}
where $\mathcal{H}_{\gamma+2r}(x;0)$ is the cdf of Chi-square
distribution with $\gamma+2r$ d.f.

Following Judge and Bock (1978) we may write the following lemma
where $E_{\N}$ stands for the expectation w.r.t normal theory.
About new insights, Lemma 2.1 in the below, gives possible
extension for evaluating specific moments under elliptical
assumption. As it can be seen from the result of Section 3, it is
always needed to obtain the expectation of a measurable function
for risk functions. To ease the understanding of what Lemma 2.1
proposes, consider a random elliptical variable
$\bX\sim\E_q(\bmu,\bV,g)$, then for a Borel measurable function
$h$, we have
\begin{eqnarray*}
E(h(\bX))=\int_0^\infty E_{\N}(h(\bX))\W(t)dt,
\end{eqnarray*}
where $E_{\N}(h(\bX))$ shows taking the expectation of $h$ under
the model $\N_q(\bmu,\bV)$.

For the proof of relevant expectations in Lemma 2.1 under
normality assumption see Judge and Bock (1978) and Saleh (2006).
\begin{lem}
If $\bw\sim \E_q(\bet,\sigma^2\bI_q,g)$, and $\Phi$ is a
measurable function then
\begin{enumerate}
  \item[(i)]
\begin{eqnarray}
E[\bw\Phi (\bw'\bw)] & = & \int_0^\infty E_{\N} [\bw\Phi (\bw'\bw)|t]\; \W(t)\;dt \nonumber \\
{} & = & \bet \int_0^\infty E_{\N} [\Phi
\left(\chi_{q+2}^2(\Delta^2_t)\right)|t]\; \W(t)\;dt,
\end{eqnarray}

  \item[(ii)]
\begin{eqnarray}
E[(\bw \bw')\Phi (\bw'\bw)] & = & \bI_q\int_0^\infty E_{\N} [\Phi
\left(\chi_{q+2}^2(\Delta^2_t)\right)|t]\; \W(t)\;dt \nonumber \\
{} & + & \bet'\bet \int_0^\infty E_{\N} [\Phi
\left(\chi_{q+4}^2(\Delta^2_t)\right)|t]\; \W(t)\;dt,
\end{eqnarray}

  \item[(iii)]
\begin{eqnarray}
E[(\bw' \bA\; \bw)\Phi (\bw'\bw)] & = &
\boldsymbol{tr}(\bA)\int_0^\infty E_{\N}
[\Phi\left(\chi_{q+2}^2(\Delta^2_t)\right)|t]\; \W(t)\;dt \nonumber \\
{} & + & (\bet \bA '\bet) \int_0^\infty
E_{\N}[\Phi\left(\chi_{q+2}^2(\Delta^2_t)\right)|t]\; \W(t)\;dt,
\end{eqnarray}
\end{enumerate}
where $\Delta^2_t=t\bet'\bet$, $\bA\in S(q)$.\\
Further,
\begin{enumerate}
  \item [(i)]
\begin{equation}
E^{(2-h)}[\chi_{q+s}^{*^{-2}}(\Delta_*^2)]=\sum_{r\geq 0}
\frac{1}{r!}K_r^{(h)}(\Delta_*^2)(q+s-2+2r)^{-1},
\end{equation}
  \item [(ii)]
\begin{equation}
E^{(2-h)}[\chi_{q+s}^{*^{-4}}(\Delta_*^2)]=\sum_{r\geq 0}
\frac{1}{r!}K_r^{(h)}(\Delta_*^2)(q+s-2+2r)^{-1}(q+s-4+2r)^{-1},
\end{equation}
\end{enumerate}
and for $h=0,1$
\begin{eqnarray}\label{eq219}
K^{(h)}_r(\Delta_*^2)&=&[-2\psi'(0)]^r\left(\frac{\Delta_*^2}{2}\right)^r\int^{\infty}_0\frac{(t^{-1})^{-r+h}}{r!}\;e^{\frac{-t\Delta_*^2
[-2\psi'(0)]}{2}}\W(t)dt,
\end{eqnarray}
If $\Phi(\bw'\bw)=I(\bw'\bw < c)$, where $I(A)$ is the indicator
function of the set $A$, then
\begin{enumerate}
  \item[(i)]
\begin{eqnarray}
E[\bw\; I(\bw'\bw<c)] & = & \bet \int_0^\infty E_{\N} [I
\left(\chi_{q+2}^2(\Delta_*^2)<c\right)|t]\; \W(t)\;dt \nonumber \\
{} & = & \bet\; \mathcal{H}_{q+2}^*(c,\Delta_*^2),
\end{eqnarray}

  \item[(ii)]
\begin{eqnarray}
E[\bw\bw'\; I(\bw'\bw<c)] & = & \bI_q \;
\mathcal{H}_{q+2}^*(c,\Delta_*^2)+\bet\bet' \;
\mathcal{H}_{q+4}^*(c,\Delta_*^2),
\end{eqnarray}

  \item[(iii)]
\begin{eqnarray}
E[\bw'\bA\; \bw\; I(\bw'\bw<c)] & = & tr(\bA) \;
\mathcal{H}_{q+2}^*(c,\Delta_*^2)\nonumber\\
{} & {}& +\;\bet'\bA\; \bet\; \;
\mathcal{H}_{q+4}^*(c,\Delta_*^2),
\end{eqnarray}
\end{enumerate}
Further
\begin{enumerate}
  \item[(i)]
\begin{eqnarray}\label{GG}
E\left[\bw\; I\left(\frac{\bw'\bw}{q S^2}<c\right)\right] & = &
\bet G_{q+2,m}^{(2)}\left(\frac{q}{q+2} c;\Delta_*^2\right),\qquad
m=n-p,
\end{eqnarray}

  \item[(ii)]
\begin{eqnarray}
E\left[\bw\bw'\; I\left(\frac{\bw'\bw}{q S^2}<c\right)\right] & =
& \bI_q
G_{q+2,m}^{(2)}\left(\frac{q}{q+2} c;\Delta_*^2\right)\nonumber\\
{} & {} & +\;\bet\bet'G_{q+4,m}^{(2)}\left(\frac{q}{q+2}
c;\Delta_*^2\right),
\end{eqnarray}

  \item[(iii)]
\begin{eqnarray}
E\left[\bw'\bA\; \bw\; I\left(\frac{w'w}{q S^2}<c\right)\right] &
= & tr(\bA)
G_{q+2,m}^{(2)}\left(\frac{q}{q+2} c;\Delta_*^2\right)\nonumber\\
{} & {} & +\;\bet' \bA\; \bet\; G_{q+4,m}^{(2)}\left(\frac{q}{q+2}
c;\Delta_*^2\right).
\end{eqnarray}
\end{enumerate}
Further,
\begin{enumerate}
  \item[(i)]
\begin{eqnarray}\label{eq46}
&&G^{(2-h)}_{q+2i,n-p}(x',\Delta_*^2)=\sum_{r=0}^{\infty}K^{(h)}_r(\Delta_*^2)I_{x'}\left[\frac{q+2i}{2}+r,\frac{n-p}{2}\right],\\
&&\hspace{.5cm}x'=\frac{qF_{\alpha}}{n-p+qF_{\alpha}},
\end{eqnarray}

   \item[(ii)]
\begin{eqnarray}\label{eq48}
&&E^{(2-h)}[F^{-j}_{q+s,n-p}(\Delta_*^2)I(F_{q+s,n-p}(\Delta_*^2)<\frac{qd}{q+s})]\\
&&=\sum_{r=0}^{\infty}K^{(h)}_r(\Delta_*^2)\bigg(\frac{q+s}{n-p}\bigg)^j\frac{B(\frac{q+s+2r-2j}{2},
\frac{n-p+2j}{2})}{B(\frac{q+s+2r}{2},\frac{n-p}{2})}I_{x}\left[\frac{q+s+2r-2j}{2},\frac{n-p+2j}{2}\right],\nonumber\\
&&x=\frac{qd}{n-p+qd}.\nonumber
\end{eqnarray}
\end{enumerate}
\end{lem}
\noindent In many practical situations, along with the model one
may suspect that $\bb$ belongs to the sub-space defined by
$\bH\bb=\bbh$. In such situation one combines the estimate of
$\bb$ and the test-statistic to obtain 3 or more estimators as in
Saleh (2006), in addition to the unrestricted and the restricted
estimators of $\bb$. Now, we proceed to define three more
estimators of $\bb$ combining the unrestricted, restricted and the
test-statistic $\nl$ as in Saleh (2006). First we consider the
preliminary test estimator (PTE) of $\bb$ which is a convex
combination of $\bt$ and $\bh$:
\begin{equation}\label{226}
\bh^{PT} = \bt I (\nl \geq F_\alpha) + \bh I (\nl < F_\alpha),
\end{equation}
where $I(A)$ is the indicator function of the set $A$ and
$F_\alpha$ is the upper $\alpha^{th}$ percentile of the central
F-distribution with $(q,m)$ d.f.

The PTE has the disadvantage that it depends on $\alpha\; (0<
\alpha < 1)$, the level of significance and also it yields the
extreme results, namely $\bh$ and $\bt$ depending on the outcome
of the test. Therefore, we define Stein-type shrinkage estimator
(SE) of $\beta$, as
\begin{eqnarray}\label{eq227}
\bh^S &  = &  \bh + (1- d\nl^{-1}) (\bt - \bh)\nonumber\\
{} & = & \bt - d\nl^{-1} (\bt - \bh),
\end{eqnarray}
where
\begin{equation}\label{eq228}
d=\frac {(q-2)(n-p)}{q(n-p+2)},\qquad \textrm{and} \qquad q \geq
3.
\end{equation}

The SE has the disadvantage that it has strange behavior for small
values of $\nl$. Also, the shrinkage factor $(1- d\nl^{-1})$ becomes
negative for $\nl < d$. Hence we define a better estimator by
positive-rule shrinkage estimator (PRSE) of $\beta$ as
\begin{eqnarray}\label{eq229}
\bh^{S+} &  = &  \bh + (1- d\nl^{-1})I[\nl > d] (\bt - \bh)\nonumber\\
{} & = & \bh^{S} - (1- d\nl^{-1})I[\nl \leq d] (\bt - \bh).
\end{eqnarray}
Note that this estimator is also a convex combination of $\bh$ and
$\bt$.

The biases, the quadratic risks and MSE-matrices of the estimators
are given in the following section and the dominance properties
are studied in section 4.

\setcounter{equation}{0}
\section{Bias and Quadratic Risk of the Estimators}
Consider for a given non-singular matrix $\bW\in S(p)$, the
\emph{weighted quadratic error loss function} of the form
\begin{equation}\label{eq31}
\boldsymbol{L}(\bb^*;\bb)=(\bb^* - \bb)'\bW(\bb^* - \bb),
\end{equation}
where $\bb^*$ is any estimator of $\bb$. Then the weighted
quadratic risk function associated with \eqref{eq31} is defined as
\begin{equation}\label{eq32}
\boldsymbol{R}(\bb^*;\bb)=E[(\bb^* - \bb)'\bW(\bb^* - \bb)].
\end{equation}

In this section, we determine the biases, and using the risk
function \eqref{eq32}, evaluate the quadratic risks and MSE
matrices of the five different estimators under study. First we
consider bias expressions of the estimators.
\newcommand{\G}{G}

Directly
\begin{equation}\label{eq33}
\boldsymbol{b}_1=E[\bt-\bb]=\boldsymbol{0}.
\end{equation}
Also
\begin{equation}\label{eq34}
\boldsymbol{b}_2=E[\bh-\bb]=-\bd.
\end{equation}
Using Lemma 2.1 we have
\begin{eqnarray}
\boldsymbol{b}_3&=&E(\bh^{PT}-\bb)\nonumber\\
 & = & E[\bt-I(\nl \leq F_\alpha)(\bt -\bh)-\bb]\nonumber\\
 & = & E[\bt-\bb]-E[I(\nl \leq F_\alpha)(\bt -\bh)]\nonumber\\
 & = & -\bC\bH'\bV_1^{1/2}\;E[I(\nl \leq F_\alpha)\bV_1^{1/2}(\bH\bt -\bbh)]\nonumber\\
 & = & -\bd \G_{q+2,m}^{(2)}\left(F_\alpha;\Delta_*^2\right).
\end{eqnarray}
Applying Lemma 2.1 we obtain
\begin{eqnarray}
\boldsymbol{b}_4 & = & E(\bh^{S}-\bb)\nonumber\\
 & = & E[\bt-d \nl^{-1}(\bt -\bh)-\bb]\nonumber\\
 & = & -d\bC^{-1}\bH'\bV_1^{1/2}E[\nl^{-1}\bV_1^{1/2}(\bH\bt-\bbh)]\nonumber\\
 & = & -dq\bd E^{(2)}[\chi_{q+2}^{*^{-2}}(\Delta_*^2)],
\end{eqnarray}
Finally, by making use of Lemma 2.1 once more, we get
\begin{eqnarray}
\boldsymbol{b}_5 & = & E(\bh^{S+}-\bb)\nonumber\\
 & = & E(\bh^{S}-\bb)-E[I(\nl\leq d)(\bt-\bh)]+d E[\nl^{-1}I(\nl\leq d)(\bt-\bh)]\nonumber\\
 & = & -dq\delta E^{(2)}[\chi_{q+2}^{*^{-4}}(\Delta_*^2)]+\delta G_{q+2,m}^{(2)}\left(d;\Delta_*^2\right)\nonumber\\
 &&+\frac{q d}{q+2}\bd E^{(2)}\left[F^{-1}_{q+2,m}(\Delta^2_\ast)I\left(F_{q+2,m}(\Delta^2_\ast\right)\leq
\frac{qd}{q+2})\right].
\end{eqnarray}
Note that as the non-centrality parameter
$\Delta_\ast^2\rightarrow\infty$,
$\boldsymbol{b}_1=\boldsymbol{b}_3=\boldsymbol{b}_4=\boldsymbol{b}_5=\0$
while $\boldsymbol{b}_2$ becomes unbounded. However, under
$H_0:\bH\bb=\bbh$, because $\bd=\0$,
$\boldsymbol{b}_1=\boldsymbol{b}_2=\boldsymbol{b}_3=\boldsymbol{b}_4=\boldsymbol{b}_5=\0$.\\

For the risks of the estimators, considering quadratic risk
function given by (3.2), we get
\begin{eqnarray}
\boldsymbol{R}(\bt;\bb)&=&E[(\bt-\bb)'\bW(\bt-\bb)]\nonumber\\
&=&\sigma^2_{\be} tr(\bC^{-1}\bW).
\end{eqnarray}
Also
\begin{eqnarray}
\boldsymbol{R}(\bh;\bb)&=&E[(\bh-\bb)'\bW(\bh-\bb)]\nonumber\\
&=&\sigma^2_{\be}tr(\bV_2\bW)+\bd'\bW\bd\nonumber\\
&=&R(\bt;\bb)-\sigma^2_{\be}tr[\bW\bC \bH'\bV_1\bH\bC]+\bd'\bW\bd.
\end{eqnarray}
Note that $\bR=\bC_1^{-1/2}\bH'\bV_1\bH\bC_1^{-1/2}$ is a
symmetric idempotent matrix of rank $q\leq p$. Thus, there exist
an orthogonal matrix $\bQ$ (see Searle, 1982) such that
\begin{eqnarray}
\bQ\bR\bQ^\prime= \left[ \begin{array}{r r}
      \bI_q & \0 \\ \0 & \0
             \end{array} \right], \ \ \bQ\bC_1^{-1/2}\bW\bC_1^{-1/2}\bQ^\prime= \left[ \begin{array}{r r}
      \bA_{11} & \bA_{12} \\ \bA_{21} & \bA_{22}
             \end{array} \right]=\bA.
\end{eqnarray}
The matrices $\bA_{11}$ and $\bA_{22}$ are of order $q$ and $p-q$
respectively.\\
Now we define random variable
\begin{eqnarray}
\bw&=&\bQ\bC_1^{1/2}\bt-\bQ\bC_1^{-1/2}\bH'\bV_1\bbh,
\end{eqnarray}
then, $\bw\sim \E_p(\bet,\sigma^2_{\be}\bI_p,g)$, where
\begin{eqnarray}
\bet&=&\bQ\bC_1^{1/2}\bb-\bQ\bC_1^{-1/2}\bH'\bV_1\bbh.
\end{eqnarray}
Partitioning the vectors $\bw=(\bw_1^\prime,\bw_2^\prime)^\prime$
and $\bet=(\bet_1^\prime,\bet_2^\prime)^\prime$ where $\bw_1$ and
$\bw_2$ are sub-vectors of order $q$ and $p-q$ respectively, we
obtain
\begin{eqnarray}
\bt-\bb&=&\bC_1^{-1/2}\bQ'(\bw-\bet).
\end{eqnarray}
Thus, we can rewrite
\begin{eqnarray}\label{eq314}
\nl=\frac{\bw_1^\prime \bw_1}{q S^2}, \quad
\bet_1'\bet_1=\theta=(\bH\bb-\bbh)'\bV_1(\bH\bb-\bbh).
\end{eqnarray}
Also we get
\begin{eqnarray}\label{eq315}
tr\{\bW[\bC \bH'\bV_1\bH\bC]\}&=&tr\{\bQ\bC_1^{-1/2}\bW\bC_1^{-1/2}\bQ'\bQ\bR\bQ'\}\nonumber\\
&=&tr\left(\left[ \begin{array}{r r}
      \bA_{11} & \bA_{12} \\ \bA_{21} & \bA_{22}
             \end{array} \right]\left[ \begin{array}{r r}
      \bI_q & \0 \\ \0 & \0
             \end{array} \right]\right)\nonumber\\
&=&tr(\bA_{11}).
\end{eqnarray}
And
\begin{eqnarray}\label{eq316}
\bd^\prime \bW\bd&=&(\bH\bb-\bbh)'\bV_1\bH\bC \bW\bC \bH'\bV_1(\bH\bb-\bbh)\nonumber\\
&=&\bet_1^\prime \bA_{11}\bet_1.
\end{eqnarray}
Substituting \eqref{eq315}-\eqref{eq316} in (3.9) yields
\begin{eqnarray}
R(\bh;\bb)&=&R(\bt;\bb)-\sigma^2_{\be}tr(\bA_{11})+\bet_1^\prime
\bA_{11}\bet_1.
\end{eqnarray}
Similarly, we obtain
\begin{eqnarray*}
R(\bh^{PT};\bb)&=&E[(\bh^{PT}-\bb)'\bW(\bh^{PT}-\bb)]\\
 & = & E\{[(\bt-\bb)-I(\nl\leq
F_{\alpha})(\bt-\bh)]'\bW[(\bt-\bb)\\
 && -I(\nl\leq F_{\alpha})(\bt-\bh)]\}\\
 & = & E[(\bt-\bb)'\bW(\bt-\bb)]-2E[I(\nl\leq
F_{\alpha})(\bt-\bb)'\bW(\bt-\bh)]\\
 &&+E[I(\nl\leq F_{\alpha})(\bt-\bh)'\bW(\bt-\bh)]\\
 &=&R(\bt;\bb)-E[\bw_1'
\bA_{11}\bw_1I(\nl\leq F_{\alpha})]\\
 && -2E[\bw_2'\bA_{21}\bw_1I(\nl\leq F_{\alpha})]+2\bet_1'
\bA_{11}E[\bw_1I(\nl\leq F_{\alpha})]\\
 && +2\bet_2'\bA_{21}E[\bw_1I(\nl\leq F_{\alpha})].
\end{eqnarray*}
By making use of the representation given by (1.4), we conclude
that $\bw_1\mid t$ and $\bw_2\mid t$ are conditionally independent
under normal theory, then one may write
\begin{eqnarray*}
E[\bw_2'\bA_{21}\bw_1I(\nl\leq F_{\alpha})]&=&\int_0^\infty
\W(t)E_N[\bw_2'\bA_{21}\bw_1I(\nl\leq F_{\alpha})] dt\\
&=&\bet_2^\prime \bA_{21}\int_0^\infty \W(t)E_N[\bw_1I(\nl\leq
F_{\alpha})] dt.
\end{eqnarray*}
Therefore, we get
\begin{eqnarray}
R(\bh^{PT};\bb)
 & = &
 R(\bt;\bb)-\sigma^2_{\be}tr(\bA_{11})\G^{(1)}_{q+2,m}\left(F_\alpha;\Delta_\ast^2\right)\nonumber\\
 && +\bet_1'\bA_{11}\bet_1\left[2\G^{(2)}_{q+2,m}\left(F_\alpha;\Delta_\ast^2\right)
 -\G^{(2)}_{q+4,m}\left(F_\alpha;\Delta_\ast^2\right)\right].~~~~~~
\end{eqnarray}
Now consider for every $q\geq3$, one can obtain
\begin{eqnarray}
E^{(2)}[\chi_{q}^{*^{-2}}(\Delta_*^2)]-E^{(2)}[\chi_{q+2}^{*^{-2}}(\Delta_*^2)]&=&2E^{(2)}[\chi_{q+2}^{*^{-4}}(\Delta_*^2)],\\
E^{(1)}[\chi_{q+2}^{*^{-2}}(\Delta_*^2)]-(q-2)E^{(1)}[\chi_{q+2}^{*^{-4}}(\Delta_*^2)]&=&\Delta^2_\ast
E^{(2)}[\chi_{q+4}^{*^{-4}}(\Delta_*^2)].
\end{eqnarray}
Then, one can find
\begin{eqnarray}
R(\bh^S;\bb)&=&E[(\bh^S-\bb)'\bW(\bh^S-\bb)]\nonumber\\
 & = &
 E[(\bt-\bb)'\bW(\bt-\bb)]-2dE[\nl^{-1}(\bt-\bb)'\bW(\bt-\bh)]\nonumber\\
 && +d^2E[\nl^{-2}(\bt-\bh)'\bW(\bt-\bh)]\nonumber\\
 & = & R(\bt;\bb)-2dE[\nl^{-1}(\bw_1'\bA_{11}\bw_1-\bet_1'\bA_{11}\bw_1+\bw_2'\bA_{21}\bw_1\nonumber\\
 && -\bet_2'\bA_{21}\bw_1)]+d^2E[\nl^{-2}(\bw_1'\bA_{11}\bw_1)]\nonumber\\
 & = & R(\bt;\bb)-dq\sigma^2_{\be}tr(\bA_{11})\bigg\{(q-2)E^{(1)}[\chi_{q+2}^{*^{-4}}(\Delta_*^2)]\nonumber\\
 &&+\left[1-\frac{(q+2)\bet_1'\bA_{11}\bet_1}{2\sigma^2_{\be}\Delta^2_\ast tr(\bA_{11})}\right]
 (2\Delta^2_\ast)E^{(2)}[\chi_{q+4}^{*^{-4}}(\Delta_*^2)]\bigg\}.
 \end{eqnarray}
Finally, the risk of PRSE is given by
\begin{eqnarray}
R(\bh^{S+};\bb)&=&E[(\bh^{S+}-\bb)'\bW(\bh^{S+}-\bb)]\nonumber\\
 & = & R(\bh^S;\bb)+E[(1-d\nl^{-1})^2I(\nl\leq
 d)(\bt-\bh)'\bW(\bt-\bh)]\nonumber\\
 && -2E[(1-d\nl^{-1})I(\nl\leq
 d)(\bh^S-\bb)'\bW(\bt-\bh)].
\end{eqnarray}
But using the fact that
\begin{eqnarray*}
 && E[\bh^S-\bb)'\bW(1-\rho\nl^{-1})I(\nl\leq d)(\bt-\bh)]\nonumber\\
 & = & E[(\bh-\bb)+(1-d\nl^{-1})(\bt-\bh)]'\bW(1-d\nl^{-1})I(\nl\leq d)(\bt-\bh)]\nonumber\\
 & = & E\{(\bh-\bb)'\bW[(1-d\nl^{-1})I(\nl\leq d)(\bt-\bh)]\}\nonumber\\
 && +E[(1-d\nl^{-1})^2I(\nl\leq d)(\bt-\bh)'\bW(\bt-\bh)],
\end{eqnarray*}
we get
\begin{eqnarray}
R(\bh^{S+};\bb)&=&R(\bh^S;\bb)-E[(1-d\nl^{-1})^2I(\nl\leq d)(\bt-\bh)'\bW(\bt-\bh)]\nonumber\\
 &&-2E\{(\bh-\bb)'\bW[(1-d\nl^{-1})I(\nl\leq d)(\bt-\bh)]\}\nonumber\\
 & = & R(\bh^S;\bb)-E[(1-d\nl^{-1})^2I(\nl\leq d )\bw_1'\bA_{11}\bw_1]\nonumber\\
 &&-2E[(1-d\nl^{-1})I(\nl\leq d)(\bw_1'\bA_{11}\bw_1-\bet_1'\bA_{11}\bw_1+\bw_2'\bA_{21}\bw_1\nonumber\\
 &&-\bet_2'\bA_{21}\bw_1)]\nonumber\\
 & = &
 R(\bh^S;\bb)\nonumber\\
 &&-\sigma^2_{\be}\bigg\{tr(\bA_{11})E^{(1)}\left[(1-\frac{q d}{q+2}F^{-1}_{q+2,m}(\Delta^2_\ast))^2
 I(F_{q+2,m}(\Delta^2_\ast)\leq
\frac{q d}{q+2})\right]\nonumber\\
 &&+\frac{\bet_1'\bA_{11}\bet_1}{\sigma^2_{\be}}E^{(2)}\left[(1-\frac{q d}{q+2}F^{-1}_{q+2,m}(\Delta^2_\ast))^2
 I(F_{q+2,m}(\Delta^2_\ast)\leq\frac{q d}{q+2})\right]\bigg\}\nonumber\\
 && -2\bet_1'\bA_{11}\bet_1E^{(2)}\left[(1-\frac{q d}{q+2}F^{-1}_{q+2,m}(\Delta^2_\ast))I(F_{q+2,m}(\Delta^2_\ast)\leq
\frac{qd}{q+2})\right].~~~~~~~~~
\end{eqnarray}
\setcounter{equation}{0}
\section{Comparison}
Providing risk analysis of the underlying estimators with the
weight matrix $\bW$, from Theorem A.2.4. of Anderson (2003), we
have
\begin{eqnarray*}
\theta ch_1(\bA_{11})\leq\bet_1^\prime \bA_{11}\bet_1\leq\theta
ch_q(\bA_{11}),
\end{eqnarray*}
or equivalently
\begin{eqnarray*}
\sigma^2_{\be}\Delta_\ast^2 ch_1(\bA_{11})\leq\bet_1^\prime
\bA_{11}\bet_1\leq \sigma^2_{\be}\Delta_\ast^2 ch_q(\bA_{11}),
\end{eqnarray*}
where $ch_1(\bA_{11})$ and $ch_q(\bA_{11})$ are the minimum and
maximum eigenvalue of $\bA_{11}$ respectively and
$\Delta_\ast^2=\theta/\sigma^2_{\be}$.

Then by (3.8) and (3.17) one may easily see that
\begin{eqnarray}
R(\bh;\bb)&\geq&
R(\bt;\bb)-\sigma^2_{\be}tr(\bA_{11})+\sigma^2_{\be}\Delta_\ast^2 ch_1(\bA_{11}),\nonumber\\
R(\bh;\bb)&\leq& R(\bt;\bb)-\sigma^2_{\be}
tr(\bA_{11})+\sigma^2_{\be}\Delta_\ast^2 ch_q(\bA_{11}).
\end{eqnarray}
When $\Delta_\ast^2$ is equal to zero, we have the above in
equalities. Thus the restricted estimator ($\bh$) dominates the
generalized least square estimator ($\bt$) denoting by
$\bh\succeq\bt$ whenever
\begin{eqnarray}
\Delta_\ast^2&\leq&\frac{tr(\bA_{11})}{ch_q(\bA_{11})},
\end{eqnarray}
while $\bh\preceq\bt$ whenever
\begin{eqnarray}
\Delta_\ast^2&\geq&\frac{tr(\bA_{11})}{ch_1(A_{11})}.
\end{eqnarray}
For $\bW=\bC$ because $tr(\bA_{11})=q$ we conclude that $\bh$
performs better than $\bt$ ($\bh\succeq\bt$) in the interval
$[0,q\sigma^2_{\be})$ and worse outside this interval.\\
Comparing $\bh^{PT}$ versus $\bt$, using risk difference, we have
\begin{eqnarray}\label{eq44}
R(\bt;\bb)-R(\bh^{PT};\bb)&=&\sigma^2_{\be}tr(\bA_{11})
G^{(1)}_{q+2,m}\left(F_\alpha;\Delta_\ast^2\right)\nonumber\\
&&-\bet_1'\bA_{11}\bet_1[2G^{(2)}_{q+2,m}\left(F_\alpha;\Delta_\ast^2\right)-G^{(2)}_{q+4,m}\left(F_\alpha;\Delta_\ast^2\right)].\nonumber\\
\end{eqnarray}
It follows that the right hand side of \eqref{eq44} is nonnegative
whenever
\begin{eqnarray}
\Delta_\ast^2&\leq&\frac{tr(\bA_{11})}{ch_q(\bA_{11})}\times\frac{G^{(1)}_{q+2,m}\left(F_\alpha;\Delta_\ast^2\right)}
{2G^{(2)}_{q+2,m}\left(F_\alpha;\Delta_\ast^2\right)-G^{(2)}_{q+4,m}\left(F_\alpha;\Delta_\ast^2\right)},
\end{eqnarray}
and vice versa. Also Under $H_0:\bH\bb=\bbh$, $\bh^{PT}\succeq\bt$.\\
Now we compare $\bh$ and $\bh^{PT}$ by the risk difference as
follows
\begin{eqnarray}
R(\bh;\bb)-R(\bh^{PT};\bb)&=&-\sigma^2_{\be}tr(\bA_{11})[1-
G^{(1)}_{q+2,m}\left(F_\alpha;\Delta_\ast^2\right)]\nonumber\\
&&+\bet_1'\bA_{11}\bet_1[1-2G^{(2)}_{q+2,m}\left(F_\alpha;\Delta_\ast^2\right)\nonumber\\
&&+G^{(2)}_{q+4,m}\left(F_\alpha;\Delta_\ast^2\right)].
\end{eqnarray}
Thus $\bh\succeq\bh^{PT}$ whenever
\begin{eqnarray}
\Delta_\ast^2&\leq&\frac{tr(\bA_{11})}{ch_q(\bA_{11})}
\times\frac{1-G^{(1)}_{q+2,m}\left(F_\alpha;\Delta_\ast^2\right)}
{1-2G^{(2)}_{q+2,m}\left(F_\alpha;\Delta_\ast^2\right)+G^{(2)}_{q+4,m}\left(F_\alpha;\Delta_\ast^2\right)},
\end{eqnarray}
and vice versa. However, under $H_0$, the dominance order of
$\bt$, $\bh$ and $\bh^{PT}$ is as follows
\begin{eqnarray}
\bh\succeq\bh^{PT}\succeq\bt.
\end{eqnarray}
In order to determine the superiority of $\bh^S$ to $\bt$, it is
enough to see that the following risk difference
\begin{eqnarray}
R(\bt;\bb)-R(\bh^S;\bb)&=&\sigma^2_{\be}dq tr(\bA_{11})\bigg\{
(q-2)E^{(1)}[\chi_{q+2}^{*^{-4}}(\Delta_*^2)]\nonumber\\
&&+\left[1-\frac{(q+2)\bet_1'\bA_{11}\bet_1}{2\sigma^2_{\be}\Delta^2_\ast
tr(\bA_{11})}\right](2\Delta^2_\ast)
E^{(2)}[\chi_{q+4}^{*^{-4}}(\Delta_*^2)],
\end{eqnarray}
is positive for all $\mathbb{A}$ such that
\begin{eqnarray}
\left\{\mathbb{A}:\frac{tr(\bA_{11})}{ch_q(\bA_{11})}\geq\frac{q+2}{2}\right\},
\end{eqnarray}
which asserts $\bh^S$ uniformly dominates $\bt$.

Further, we show that the shrinkage factor $d$ of the Stein-type
estimator is robust with respect to $\bb$ and the unknown mixing
distribution.
\begin{thm}
Consider the model \eqref{eq11} where the error-vector belongs to
the ECD, $\E_n(\boldsymbol{0},\sigma^2\bV,g)$. Then the Stein-type
shrinkage estimator, $\bh^S$ of $\bb$ given by
\begin{eqnarray*}
\bh^S=\bt - d^*\nl^{-1} (\bt - \bh)
\end{eqnarray*}
uniformly dominates the unrestricted estimator $\bt$ with respect
to the quadratic loss function given by \eqref{eq31} for
$\bW=\bC$, and is minimax if and only if $0<
d^*\leq\frac{2m}{m+2}$. The largest reduction of the risk is
attained when $d^*=\frac{m}{m+2}$.
\end{thm}
For the proof see the Appendix.
\begin{rem}
Consider the coefficient $d$ given by (2.31). From $q\geq3$, we
get $0<d=\frac{(q-2)m}{q(m+2)}<\frac{2m}{m+2}$ and thus using
Theorem 4.1, $\bh^S$ in equation (2.30) uniformly dominates $\bt$
on the whole parameter space under quadratic loss function.
\end{rem}
To compare $\bh$ and $\bh^S$ we may write
\begin{eqnarray}
R(\bh^S;\bb)&=&R(\bh;\bb)+\sigma^2_{\be}tr(\bA_{11})-\bet_1'\bA_{11}\bet_1-dq\sigma^2_{\be}tr(\bA_{11})\bigg\{
(q-2)E^{(1)}[\chi_{q+2}^{*^{-4}}(\Delta_*^2)]\nonumber\\
&&+\left[1-\frac{(q+2)\bet_1'\bA_{11}\bet_1}{2\sigma^2_{\be}\Delta^2_\ast
tr(\bA_{11})}\right](2\Delta^2_\ast)
E^{(2)}[\chi_{q+4}^{*^{-4}}(\Delta_*^2)].
\end{eqnarray}
Under $H_0$, this becomes
\begin{eqnarray}
R(\bh^S;\bb)&=&R(\bh;\bb)+\sigma^2_{\be}(1-d)tr(\bA_{11})\nonumber\\
&&\geq R(\bh;\bb),
\end{eqnarray}
while
\begin{eqnarray}
R(\bh;\bb)&=&R(\bt;\bb)-\sigma^2_{\be}tr(\bA_{11})\nonumber\\
&&\leq R(\bt;\bb).
\end{eqnarray}
Therefore, $\bh$ performs better that $\bh^S$ under $H_0$.
However, as $\bet_1$ moves away from $\0$, $\bet_1'\bA_{11}\bet_1$
increases and the risk of $\bh$ becomes unbounded while the risk
of $\bt^S$ remains below the risk of $\bt$; thus $\bt^S$ dominates
$\bh$ outside an interval around the origin.\\
Consider under $H_0$
\begin{eqnarray}
R(\bh^S;\bb)&=&R(\bh^{PT};\bb)+\sigma^2_{\be}tr(\bA_{11})[1-\alpha-d]\nonumber\\
&&\geq R(\bh^{PT};\bb),
\end{eqnarray}
for all $\alpha$ such that
$F^{-1}_{q+2,m}(d,0)\leq\frac{qF_\alpha}{q+2}$. This means the
estimator $\bh^S$ does not always dominates $\bh^{PT}$ under
$H_0$.\\
Thus, under $H_0$ with $\alpha$ satisfying
$F^{-1}_{q+2,m}(d,0)\leq\frac{qF_\alpha}{q+2}$, we can conclude
that
\begin{eqnarray}
\bh\succeq\bh^{PT}\succeq\bh^S\succeq\bt.
\end{eqnarray}

We compare, the risks of $\bh^{S+}$ and $\bh^S$. The risk
difference is given by
\begin{eqnarray*}
R(\bh^{S+};\bb)-R(\bh^S;\bb)&=&-\sigma^2_{\be}\bigg\{tr(\bA_{11})E^{(1)}\left[(1-\frac{q
d}{q+2}F^{-1}_{q+2,m}(\Delta^2_\ast))^2
I(F_{q+2,m}(\Delta^2_\ast)\leq
\frac{q d}{q+2})\right]\nonumber\\
 &&+\frac{\bet_1'\bA_{11}\bet_1}{\sigma^2_{\be}}E^{(2)}\left[(1-\frac{q d}{q+2}F^{-1}_{q+2,m}(\Delta^2_\ast))^2
 I(F_{q+2,m}(\Delta^2_\ast)\leq\frac{q d}{q+2})\right]\bigg\}\nonumber\\
 && -2\bet_1'\bA_{11}\bet_1E^{(2)}\left[(1-\frac{q d}{q+2}F^{-1}_{q+2,m}(\Delta^2_\ast))I(F_{q+2,m}(\Delta^2_\ast)\leq
\frac{qd}{q+2})\right]
\end{eqnarray*}
The right hand side of the above equality is negative since for
$F_{q+2,m}(\Delta^2_\ast)\leq \frac{q d}{q+2}$, $(\frac{q
d}{q+2}F_{q+2,m}(\Delta^2_\ast)-1)\geq0$ and also the expectation
of a positive random variable is positive. That for all $\bb$,
$R(\bh^{S+};\bb)\leq R(\bh^{S};\bb)$. And using Theorem 4.1,
$R(\bh^{S+};\bb)\leq R(\bt;\bb)$.\\
To compare $\bh$ and $\bh^{S+}$, first consider the case under
$H_0$ i.e., $\bet_1=\0$. In this case
\begin{eqnarray}
R(\bh^{S+};\bb)&=&R(\bh;\bb)+\sigma^2_{\be} tr(\bA_{11})\bigg\{
(1-d)\nonumber\\
&&-E^{(1)}\left[(1-\frac{q d}{q+2}F^{-1}_{q+2,m}(0))^2I(F_{q+2,m}(0)\leq\frac{q d}{q+2})\right]\bigg\}\nonumber\\
&&\geq R(\bh;\bb),
\end{eqnarray}
since
\begin{eqnarray}
&&E^{(1)}\left[(1-\frac{q
d}{q+2}F^{-1}_{q+2,m}(0))^2I(F_{q+2,m}(0)\leq\frac{q
d}{q+2})\right]\nonumber\\
&&\leq E^{(1)}\left[(1-\frac{q
d}{q+2}F^{-1}_{q+2,m}(0))^2\right]=1-d.
\end{eqnarray}
Thus under $H_0$, $\bh\succeq\bh^{S+}$. However, as $\bet_1$ moves
away from $\0$, $\bet_1'\bA_{11}\bet_1$ increases and the risk of
$\bh$ becomes unbounded while the risk of $\bt^{S+}$ remains below
the risk of $\bt$; thus $\bt^{S+}$ dominates
$\bh$ outside an interval around the origin.\\
Now, we compare $\bh^{S+}$ and $\bh^{PT}$. When $H_0$ holds,
because $G^\ast_{q+2,m}(F_\alpha,0)=1-\alpha$,
\begin{eqnarray}
R(\bh^{S+};\bb)&=&R(\bh^{PT};\bb)+\sigma^2_{\be}tr(\bA_{11})\bigg\{
1-\alpha-d\nonumber\\
&&-E^{(1)}\left[(1-\frac{q
d}{q+2}F^{-1}_{q+2,m}(0))^2I(F_{q+2,m}(0)\leq\frac{q
d}{q+2})\right]\bigg\}\nonumber\\
&&\geq R(\bh^{PT};\bb),
\end{eqnarray}
for all $\alpha$ satisfying
\begin{eqnarray}
E^{(1)}\left[(1-\frac{q
d}{q+2}F^{-1}_{q+2,m}(0))^2I(F_{q+2,m}(0)\leq\frac{q
d}{q+2})\right]&\leq&1-\alpha-d.
\end{eqnarray}
Thus, $\bh^{S+}$ does not always dominates $\bh^{PT}$ when the
null-hypothesis $H_0$ holds.\\

Since always $\bh^{S+}\succeq\bh^S\succeq\bt$, and under $H_0$,
the restricted estimator $\bh$ performs better that all others,
the dominance order of the five estimators under the null
hypothesis $H_0$, can be determined under the following two
categories
\begin{eqnarray*}
1.&&\bh\succeq\bh^{PT}\succeq\bh^{S+}\succeq\bh^S\succeq\bt,\quad\mbox{satisfying}\;(4.19), \\
2.&&\bh\succeq\bh^{S+}\succeq\bh^S\succeq\bh^{PT}\succeq\bt,\quad\forall\;\alpha\;\ni\;F^{-1}_{q+2,m}(d,0)>\frac{qF_\alpha}{q+2}.
\end{eqnarray*}

\setcounter{equation}{0}
\section{Concluding Remarks}
In this paper, we proposed five different estimators for the
regression parameters of a linear regression model. In this
approach, the prior non-sample information $\bH\bb=\bbh$ is
suspected. Based on the two ordering superiority categories of the
estimators, under the constraint $\bH\bb=\bbh$, two non-linear
estimators $\bh^{S+}$ and $\bh^S$ perform better than unbiased
GLSE.

The behavior of Stein-type estimators are restricted by the
condition $q\geq3$. However, The PTE requires the size of testing
$H_0:\bH\bb=\bbh$. For $\W(t)$ as dirac delta function, the
maximal savings in risk for the shrinkage estimator is
$\frac{m(q-2)}{p(m+2)}$, while for $\W(t)$ as inverse-gamma
function, it is equal to $\frac{m(q-2)}{p(m+2)}\frac{\nu-2}{\nu}$,
where $\nu$ is d.f. of multivariate student's t-distribution.

Remarkably, the behavior of all estimators comparing with each
others under elliptical symmetry are exactly the same as under
normal theory as exhibited in Saleh (2006). This phenomenon shows
the dominance order of estimators and regarding substantial
conditions under normal theory are significantly robust.

Another noteworthy fact is that under the subclass (c) of
elliptical models we could present fundamental Theorem 4.1 under
elliptical models rather than scale mixture of normal
distributions. It is important to point out that just under the
signed measure $\W(t)$ under subclass (c) we are able to prove
Theorem 4.1 even for non-positive measures. See Srivastava and
Bilodeau (1989) for more discussion in this regard.

Finally, if $\nu_1$ remains constant, the distribution of
$F_{\nu_1,\nu_2}$ tends to that of $\chi^2_{\nu_1}/\nu_1$ as
$\nu_2$ tends to infinity (see Johnson and Kotz, 1970).\\
Now consider the class of local alternatives $\{K_{(n)}\}$
defined by
\begin{equation}
K_{(n)}: \bH\bb=\bbh+n^{\frac{-1}{2}}\boldsymbol{\xi}.
\end{equation}
Furthermore, following Saleh (2006), consider the following
regularity conditions hold\\
$(i) \max_{1\leq i\leq
n}\bx_i^\prime(\bX'\bV^{-1}\bX)^{-1}\bx_i\rightarrow0 \ \mbox{as}
\ n\rightarrow\infty$\\
where $\bx_i^\prime$ is the ith row of $\bX$;\\
$(ii) \lim_{n\rightarrow\infty} \{n^{-1}(\bX'\bV^{-1}\bX)\}=\bC$\\
for finite $\bC\in S(p)$.\\
Then using Theorem 7.8.3 from Saleh (2006) in addition to Theorem
2.1 we obtain the following important result for the test
statistic
\begin{equation}
\lim_{n\rightarrow\infty} P(\mathcal{L}_n\leq
x)=\mathcal{H}_q^*(x;\delta^2),\nonumber
\end{equation}
Based on the above results, one can easily obtain the asymptotic
distributional bias, risk and MSE matrix of each estimator under
study using the following definition
\begin{equation}
G_p(\bx)=\lim_{n\rightarrow\infty}P_{K_{(n)}}\{\sqrt{n}S^{-2}(\bb^*-\bb)\leq
\bx\}.\nonumber
\end{equation}
Then\\
$\boldsymbol{b}(\bb^*)=\int \bx dG_p(\bx)$, \
$\boldsymbol{M}(\bb^*)=\int \bx\bx'dG_p(\bx)$, \
$R(\bb^*;\bb)=tr[\bW\boldsymbol{M}(\bb^*)]$,\\
which have similar notations to those are given in this paper.
\setcounter{equation}{0}
\section{Appendix}
\textbf{Proof of Theorem 2.1}\\
The likelihood ratio is given by
\begin{eqnarray*}
\lambda&=&\frac{max_{\omega}L(\by)}{max_{\Omega}L(\by)}\\
&=&\frac{d_{n}|\hat{\sigma}^2\bV|^{-1/2}max_y
g\left[\frac{(\by-\bX\bh)'\bV^{-1}(\by-\bX\bh)}{2\hat{\sigma}^2}\right]}{d_{n}|\tilde{\sigma}^2\bV|^{-1/2}
max_y g\left[\frac{(\by-\bX\bt)'\bV^{-1}(\by-\bX\bt)}{2\tilde{\sigma}^2}\right]}\\
&=&\bigg(\frac{\tilde{\sigma}^2}{\hat{\sigma}^2}\bigg)^{n/2}\frac{g(y_g)}{g(y_g)}\\
&=&\left(\frac{1}{1+\frac{(\bH\bt-\bbh)'(\bH\bC^{-1}\bH')^{-1}(\bH\bt-\bbh)}{(\by-\bX\bt)'\bV^{-1}(\by-\bX\bt)}}\right)^{n/2}.
\end{eqnarray*}
Therefore
\begin{eqnarray*}
\lambda^{2/n}&=&\bigg(\frac{1}{1+\frac{q}{n-p}\nl}\bigg),
\end{eqnarray*}
which is decreasing with respect to $\nl$. Now consider that
\begin{eqnarray*}
(n-p)S^2|t&=&(\by-\bX\bt)'\bV^{-1}(\by-\bX\bt)|t\\
&=&\by'\bigg[\bigg(\bV^{-1}-\bV^{-1}\bX(\bX'\bV^{-1}\bX)^{-1}\bX'\bV^{-1}\bigg)\bigg]\by|t\sim\chi^2_{n-p}
\end{eqnarray*}
also $(\bH\bC^{-1}\bH')^{-1/2}(\bH\bt)|t\sim
\N_q((\bH\bC^{-1}\bH')^{-1/2}(\bH\bb-\bbh),t^{-1}\sigma^2\bI_q)$.
Then
\begin{eqnarray*}
\bt'\bH'(\bH\bC^{-1}\bH')^{-1}\bH\bt|t\sim\chi^2_{q,\Delta^2_t},
\end{eqnarray*}
where $\Delta^2_t=\frac{t\theta}{\sigma^2}$. Also using the fact
that $(\by-\bX\bt)'\bV^{-1}(\by-\bX\bt)|t$ and
$\bt'\bH'(\bH\bC^{-1}\bH')^{-1}\bH\bt|t$ are independent, we get
\begin{eqnarray*}
\nl|t=\frac{\bt'\bH'(\bH\bC^{-1}\bH')^{-1}\bH\bt|t}{qS^2|t}\sim
F_{q,n-p,\Delta^2_t}.
\end{eqnarray*}
Hence
\begin{eqnarray*}
\boldsymbol{g}^*_{q,m}(\nl)&=&\int^{\infty}_0\W(t)F_{q,m,\Delta^2_t}(\nl|t)dt\\
&=&\int^{\infty}_0\W(t)\sum^{\infty}_{r=0}\frac{e^{\frac{-\Delta^2_t}{2}}(\frac{\Delta^2_t}{2})^r}
{\Gamma(r+1)}\bigg(\frac{q}{m}\bigg)^{\frac{q}{2}+r}\nonumber\\
&&\times\frac{\nl^{\frac{q}{2}+r-1}}{B(\frac{q}{2}+r,\frac{m}{2})(1+\frac{q}{m}\nl)^{\frac{q+m}{2}+r}}\;dt\\
&=&\sum^{\infty}_{r=0}\bigg(\frac{q}{m}\bigg)^{\frac{q}{2}+r}\frac{\nl^{\frac{q}{2}+r-1}K^{(0)}_r(\Delta^2_\ast)}
{B(\frac{q}{2}+r,\frac{m}{2})(1+\frac{q}{m}\nl)^{\frac{q+m}{2}+r}}.
\end{eqnarray*}
\textbf{Proof of Theorem 4.1}\\
By making use of $\dot{\bz}=\bH'\bV_1(\bH\bt-\boldsymbol{h})$, the
SE can be rewritten as
\begin{eqnarray*}
\bh^S&=&\bt-qd^*S^2\left[(\bH\bt-\boldsymbol{h})'\bV_1(\bH\bt-\boldsymbol{h})\right]^{-1}
\bC^{-1}\bH'\bV_1(\bH\bt-\boldsymbol{h})\\
&=&\bt-qd^*S^2\left(\dot{\bz}'\bC^{-1}\dot{\bz}\right)^{-1}\bC^{-1}\dot{\bz}.
\end{eqnarray*}
Then, the risk difference of the SE and the UE under quadratic
loss function, is given by
\begin{eqnarray*}
\D_4&=&E(\bh^S-\boldsymbol{\beta})'\bC(\bh^S-\boldsymbol{\beta})-E(\bt-\boldsymbol{\beta})'\bC(\bt-\boldsymbol{\beta})\\
&=&(d^*)^2E\left[q^2S^4\left(\dot{\bz}'\bC^{-1}\dot{\bz}\right)^{-1}\right]-2d^*E\left[qS^2\left(\dot{\bz}'\bC^{-1}\dot{\bz}\right)
^{-1}(\bt-\bb)'\dot{\bz}\right]\\
&=&(d^*)^2E_{t}\left\{E_N\left[q^2S^4\left(\dot{\bz}'\bC^{-1}\dot{\bz}\right)^{-1}\bigg|t\right]\right\}\\
&&-2d^*E_{t}\left\{E_N\left[qS^2\left(\dot{\bz}'\bC^{-1}\dot{\bz}\right)
^{-1}(\bt-\bb)'\bH'\bV_1(\bH\bt-\boldsymbol{h})\bigg|t\right]\right\}\\
&=&\frac{q^2(m+2)}{m}\left(d^*\right)^2E_{t}\left(\frac{\tau^{-2}}{\dot{\bz}'\bC^{-1}\dot{\bz}}\right)
-2q^2d^*E_{t}\left(\frac{t^{-2}}{\dot{\bz}'\bC^{-1}\dot{\bz}}\right),
\end{eqnarray*}
since $\left(\frac{mS^2}{\sigma^2}\right)\bigg|t\sim
t^{-1}\chi^2_m$ and $\bt'\bH'\bV_1\bH\bt\mid t\sim
t^{-2}\sigma^4\chi^2_q(\dot{\bd})$, where
$\dot{\bd}=\bb'\bH'\bV_1\bH\bb$. \\ Therefore, $\D_4\leq0$ if and
only if $0< d^*\leq\frac{2m}{m+2}$ since $\int_0^\infty
\frac{t^{-2}}{\dot{\bz}'\bC^{-1}\dot{\bz}}\;d\W(t)>0$. (See
Srivastava and Bilodeau, 1989).
\section*{Acknowledgments}
The authors would like to thank an anonymous referee for his/her
really helpful suggestions which let to putting many details in
the paper and pointing out some mistakes in the earlier version of
the manuscript. This research was supported by a grant from
Ferdowsi University of Mashhad; (No. MS89184TAB).

\section*{References}

\baselineskip=12pt
\def\ref{\noindent\hangindent 25pt}

\ref Anderson, T. W. (2003). {\em An introduction to multivariate
statistical analysis}, 3rd ed., John Wiley and Sons, New York.

\ref Anderson, T. W., Fang, K. T. and Hsu, H. (1986).
Maximum-likelihood estimates and likelihood-ratio criteria for
multivariate elliptically contoured distributions. {\em The
Canadian J. Statist.} \textbf{14}, 55--59.

\ref Arashi, M. and Tabatabaey, S. M. M., (2010). A note on
classical Stein-type estimators in elliptically contoured models,
{\em J. Statist. Plann. Inf.}, {\bf140}, 1206-1213.

\ref Arashi, M. and Tabatabaey, S. M. M., (2008). Stein-type
improvement under stochastic constraints: Use of multivariate
Student-t model in regression, {\em Statist. Prob. Lett.},
{\bf78}, 2142-2153.

\ref Cambanis, S., Huang, S. and Simons, G. (1981). On the theory
of elliptically contoured distributions, {\em J. Mult. Anal.},
{\bf11}, 368–385.

\ref Chu, K. C., (1973). Estimation and decision for linear
systems with elliptically random process. {\em IEEE Trans. Autom.
Cont.}, {\bf18}, 499-505.

\ref Debnath, L. and Bhatta, D., (2007). {\em Integral Transforms
and Their Applications}, Chapman and Hall, London, New York.

\ref Fang, K. T., Kotz, S. and Ng, K. W. (1990). {\em Symmetric
Multivariate and Related Distributions}, Chapman and Hall, London,
New York.

\ref Gupta, A. K. and Varga, T. (1993), {\em Elliptically
Contoured Models in Statistics}, Kluwer Academic Press.

\ref Johnson, N. L. and Kotz, S. (1970), {\em Continuous
univariate distributions-2}, John Wiley and Sons, New York.

\ref Judge, G. G. and bock, M. E., (1978). {\em The statistical
implication of pre-test and Stein-rule estimators in
Econometrics}, North-Holland, New York.

\ref Kelker, D. (1970), Ditribution theory of spherical
distributions and location-scale parameter generalization, {\em
Sankhya}, {\bf 32}, 419-430

\ref Muirhead, R. J. (1982), {em Aspect of Multivariate
Statistical Theory}, John Wiley, New York.

\ref Saleh, A. K. Md. E. (2006), {\em Theory of Preliminary Test
and Stein-type Estimation with Applications}, John Wiley, New
York.

\ref Searle, S. R. (1982), {\em Matrix Algebra Useful for
Statistics}, John Wiley, New York.

\ref Srivastava, M. and Bilodeau, M., (1989). Stein estimation
under elliptical distribution, {\em J. Mult. Annal.}, {\bf28},
247-259.

\ref Tabatabaey, S. M. M., (1995), {\em Preliminary Test Approach
Estimation: Regression Model with Spherically Symmetric Errors},
Unpublished Ph.D. Thesis, Carleton University, Canada.

\end{document}